\documentclass[leqno,
12pt]{amsart}
\usepackage{amssymb}
\usepackage{ifthen}

\nonstopmode \numberwithin{equation}{section}
\newtheorem{thm}{Theorem}[section]
\newtheorem{cor}{Corollary}[section]
\newtheorem{lem}{Lemma}[section]

\newtheorem{conj}{Conjecture}

\theoremstyle{definition}

\newtheorem{prob}[equation]{Problem}

\newenvironment{rem}{%
\bigskip
\noindent \textsl{{\sl Remark. }}}{\bigskip}
\newenvironment{rems}{%
\bigskip
\noindent \textsl{{\sl Remarks. }}}{\bigskip}

\newcounter {own}
\def\theown {\thesection       .\arabic{own}}

\newenvironment{nonsec}{\bf
\setcounter{own}{\value{equation}}\addtocounter{equation}{1}
\refstepcounter{own}
\bigskip

\indent \theown.$\ \,$}{$\,$.\ \ \ }

\newenvironment{pf}[1][]{%
 \vskip 3mm
 \noindent
 \ifthenelse{\equal{#1}{}}%
  {{\slshape Proof. }}%
  {{\slshape #1.} }%
 }%
{\qed\bigskip}

\newcounter{alphabet}
\newcounter{tmp}

\newcommand{\ID}{{\mathbb D}}

\newcommand{\IC}{{\mathbb C}}

\def\be{\begin{equation}}
\def\ee{\end{equation}}

\newcommand{\bee}{\begin{enumerate}}
\newcommand{\eee}{\end{enumerate}}

\newcommand{\blem}{\begin{lem}}
\newcommand{\elem}{\end{lem}}
\newcommand{\bthm}{\begin{thm}}
\newcommand{\ethm}{\end{thm}}
\newcommand{\bcor}{\begin{cor}}
\newcommand{\ecor}{\end{cor}}
\newcommand{\beg}{\begin{examp}}
\newcommand{\eeg}{\end{examp}}
\newcommand{\begs}{\begin{examples}}
\newcommand{\eegs}{\end{examples}}
\newcommand{\bdefe}{\begin{defin}}
\newcommand{\edefe}{\end{defin}}
\newcommand{\bprob}{\begin{prob}}
\newcommand{\eprob}{\end{prob}}
\newcommand{\bei}{\begin{itemize}}
\newcommand{\eei}{\end{itemize}}

\newcommand{\bcon}{\begin{conj}}
\newcommand{\econ}{\end{conj}}
\newcommand{\bcons}{\begin{conjs}}
\newcommand{\econs}{\end{conjs}}
\newcommand{\bprop}{\begin{propo}}
\newcommand{\eprop}{\end{propo}}
\newcommand{\br}{\begin{rem}}
\newcommand{\er}{\end{rem}}
\newcommand{\brs}{\begin{rems}}
\newcommand{\ers}{\end{rems}}
\newcommand{\bo}{\begin{obser}}
\newcommand{\eo}{\end{obser}}
\newcommand{\bos}{\begin{obsers}}
\newcommand{\eos}{\end{obsers}}
\newcommand{\bpf}{\begin{pf}}
\newcommand{\epf}{\end{pf}}
\newcommand{\ba}{\begin{array}}
\newcommand{\ea}{\end{array}}
\newcommand{\beq}{\begin{eqnarray}}
\newcommand{\beqq}{\begin{eqnarray*}}
\newcommand{\eeq}{\end{eqnarray}}
\newcommand{\eeqq}{\end{eqnarray*}}

\newcommand{\ds}{\displaystyle}

\begin{document}
\bibliographystyle{amsplain}
\title{COEFFICIENT INEQUALITIES FOR
CONCAVE AND MEROMORPHICALLY STARLIKE UNIVALENT FUNCTIONS}
\author{B. Bhowmik}
\address{B. Bhowmik, Department of Mathematics,
Indian Institute of Technology Madras, Chennai-600 036, India.}
\email{ditya@iitm.ac.in}
\author{S. Ponnusamy}
\address{S. Ponnusamy, Department of Mathematics,
Indian Institute of Technology Madras, Chennai-600 036, India.}
\email{samy@iitm.ac.in}

\date{
Aug. 9, 2007
; File: bhow${}_{-}$p1${}_{-}$finalnew.tex}

\begin{abstract}
Let $\ID$ denote the open unit disk and $f:\,\ID\rightarrow\overline\IC$ be meromorphic
and univalent in $\ID$ with the simple pole at $p\in (0,1)$
and satisfying the standard normalization $f(0)=f'(0)-1=0$. Also, let $f$ have the
expansion
$$f(z)=\sum_{n=-1}^{\infty}a_n(z-p)^n,\quad |z-p|<1-p,
$$
such that $f$ maps $\ID$ onto a domain whose complement with respect to
$\overline{\IC}$ is a convex set (starlike set with respect to a point
$w_0\in \IC, w_0\neq 0$ resp.).
We call these functions as concave (meromorphically starlike resp.)
univalent functions and denote this class by $Co(p)$ $(\Sigma^s(p, w_0)$ resp.).
We prove some coefficient estimates for functions in the classes where
the sharpness of these estimates is also achieved.
\end{abstract}
\baselineskip=16pt
\thanks{}

\maketitle
\pagestyle{myheadings}
\markboth{B. Bhowmik and S. Ponnusamy}{Concave univalent function}

\section{Introduction}\label{sec1-bpw2}

One of the most interesting questions in the theory of univalent
functions is to address the region of variability of the $n$-th
Taylor (Laurent resp.) coefficient for functions $f$
that are analytic (meromorphic resp.) and univalent in the
unit disk $\ID =\{z:\,|z|<1\}$. The leading example is the
Bieberbach conjecture settled by de Branges in 1985 for the class
$\mathcal S$ of normalized analytic univalent functions $f$ in $\ID$ although
corresponding results for important subclasses of $\mathcal S$ are
relatively easy and were settled positively much earlier. In this paper,
we consider the family $Co(p)$ of
 functions $f:\ID\to \overline{\IC}$ that satisfy the following conditions:
\bee
\item[(i)] $f$ is meromorphic in $\ID$ and has a simple pole at the point
$p\in (0,1)$ with the standard normalization $f(0)=f'(0)-1=0$.
\item[(ii)] $f$ maps $\ID$ conformally onto a set whose complement with respect to
$\overline{\IC}$ is convex.
\eee
Each $f\in Co(p)$ has the power series expansion of the form
\be\label{p1eq0}
f(z)=z+\sum_{n=2}^{\infty} A_n(f)z^n, \quad |z|<p.
\ee
For our investigation, we  consider the Laurent expansion of $f\in Co(p)$
about the pole $z=p$:
\be\label{p1eq1}
f(z)=\sum_{n=-1}^{\infty} a_n(z-p)^n, \quad z\in \Delta_p,
\ee
where $\Delta_p=\{z\in \IC:\, |z-p|<1-p\}$.
Motivated by the works of Pfaltzgraff  and Pinchuk \cite{Pfa-Pin-90},
Miller \cite{mill-80}, and Livingston \cite{Living-94}, the class $Co(p)$ has been
investigated recently in \cite{A-wirths06,Avk-Wir-93,Avk-Wir-92,Wir-91,Wirths-pre}.
A necessary and sufficient condition
for a function $f$ to be in $Co(p)$  (\cite{Living-94}) is
that ${\rm Re}\,\phi(z,f) >0$ for all $z\in \ID$, where
$$
\phi(z,f)=-(1+p^2)+2pz-\frac{(z-p)(1-pz)f''(z)}{f'(z)}, \quad z\in \ID.
$$
Livingston \cite{Living-94} determined some estimates regarding the real part
of $A_n(f)$ for $n=2,3$ when $f\in Co(p)$ has
the expansion (\ref{p1eq0}). In the same article he conjectured an
estimate for the real part of the general coefficient $A_n(f)$ $(n\geq 2$)
for $ f\in Co(p)$. After a long gap of ten years,  positive developments have
occurred in this line of work. For example, the recent work of
Avkhadiev and Wirths \cite{A-wirths06} settles
the conjecture of Avkhadiev, Pommerenke and Wirths
\cite{Avk-Wir-93} which, in particular, provides a proof of the Livingston conjecture.
For a ready reference, we now recall it here.

\bigskip
\noindent
{\bf Theorem A.} {\rm \cite{A-wirths06}}
{\em Let $n\geq 2$ and $p\in (0,1).$ For each $f \in Co(p)$ with the expansion
$(1.1)$ the inequality
\be\label{p1eq21}
\left|A_n(f)-\frac{1-p^{2n+2}}{p^{n-1}(1-p^4)}\right|\leq
\frac{p^2(1-p^{2n+2})}{p^{n-1}(1-p^4)}
\ee
is valid. Equality in {\rm (\ref{p1eq21})} is attained if and only if $f$
is one of the functions $f_{\theta}$, $\theta \in [0, 2\pi)$, where
\be\label{p1eq18aa}
f_{\theta}(z)=\frac{z-\frac{p}{1+p^2}(1+e^{i\theta})z^2}{(1-\frac{z}{p})(1-zp)}.
\ee
For each complex number in the disk described in {\rm (\ref{p1eq21})} there
exists a function $f\in Co(p)$ such that this number occurs as the
$n$-th Taylor coefficient of $f.$
}
\bigskip

Interestingly,  Wirths  \cite{Wirths-pre} established the following
representation formula for functions in $Co(p)$.

\bigskip
\noindent
{\bf Theorem B.} {\rm \cite{Wirths-pre}}
{\em
For each $f\in Co(p)$, there exists a function $\omega$
holomorphic in $\ID$ such that $\omega(\ID)\subset \overline{\ID}$ and
\be\label{p1eq17}
f(z)=\frac{z-\frac{p}{1+p^2}(1+\omega(z))z^2}{(1-\frac{z}{p})(1-zp)},
\quad \mbox{$ z\in \ID $}.
\ee
}
\bigskip

The above representation formula
has been used by the authors in \cite{BPW1} to obtain some other
kind of coefficient estimates for functions in the class $Co(p)$
with the Laurent expansion of the form (\ref{p1eq1}).

In the present article, we first obtain certain coefficient estimates
for functions in $Co(p)$  but this time for the expansion of the form (\ref{p1eq1}).
Next we move on to discuss a related class of meromorphically
starlike functions, namely, the class $\Sigma^s(p, w_0$ and obtain
a simple and easily applicable
representation formula for this class. Using this formula, we also
obtain some sharp coefficient estimates for functions in this class.
As a consequence of our investigation, we rectify a mistake that appeared in
the work of Livingston in \cite[Theorem 9]{Living-94}.

Now, we state our first result.

\bthm\label{1th1}
Let $p\in(0,\frac{\sqrt 5-1}{2}]$ and $f \in Co(p)$ have the expansion $(\ref{p1eq1})$.
Then
\be\label{p1eq2}
\left |p-(1-p^2)\frac{a_0}{a_{-1}}
\right |\leq \frac{p}{|a_{-1}|} ,
\quad \mbox{i.e. }\quad
\left |a_{-1}-\frac{1-p^2}{p}a_0
\right |\leq 1.
\ee
The inequality is sharp.
\ethm

\br
In \cite{Wirths-pre} Wirths  has obtained the region of variability
for $a_{-1}(f)$, namely, the inequality
$$
\left |a_{-1}+ \frac{p^2}{1-p^4}\right |\leq \frac{p^4}{1-p^4}
\quad \mbox{for $0<p<1$}.
$$
In \cite{BPW1}, the domain of variability of $a_0(f)$
is determined by the inequality
$$\left|\frac{(1-p^2)a_0}{p}+\frac{1-p^2+p^4}{1-p^4}\right|
\leq \frac{p^2(2-p^2)}{1-p^4} \quad \mbox{for $p\in (0,\sqrt{3}-1]$}.
$$
Equality in each of the above two inequalities is attained if and only if
$f$ is one of the functions given in (\ref{p1eq18aa}).
\hfill $\Box$\er

Next result presents sharp coefficient estimates for all $n\geq 3$
if $f\in Co(p)$ has the expansion (\ref{p1eq1}).

\bthm\label{1th2}
If $f\in Co(p)$ with $p\in (0,1)$ and has the expansion $(\ref{p1eq1})$, then we have for $(n\geq 3)$\\

$\ds \left|a_{n-2}-\frac{(1-p^2)a_{n-1}}{p}\right|$
\begin{equation}\label{p1eq18}
\leq
\frac{p}{(1-p^4)(1-p)^{n-1}}
\left[1-\left (\frac{1-p^4}{p^4}\right )^2\,
\left|a_{-1}+\frac{p^2}{1-p^4}\right|^2\right].
\end{equation}
The equality holds in the above inequality for the functions $f_\theta
~~(0\leq\theta\leq 2\pi)$ of the form $(\ref{p1eq18aa})$.
\ethm


\section{Proofs of Theorems \ref{1th1} and \ref{1th2}}\label{1proofs}
\begin{nonsec}Proof of Theorem \ref{1th1}\end{nonsec}
Let $f\in Co(p)$. Then, by Theorem B,
there exists a function $\omega$ holomorphic in $\ID$ such that
$\omega(\ID)\subset \overline{\ID}$ satisfying
the representation formula (\ref{p1eq17}).

Now, let $f\in Co(p)$ have the Laurent expansion (\ref{p1eq1}) and let
$\omega$ have the Taylor expansion
\be\label{p1eq4a}
\omega(z)=\sum_{n=0}^{\infty}c_n(z-p)^n, \quad z\in\Delta_p.
\ee
Using these two expansions, the series formulation of (\ref{p1eq17}) takes
the form
\beq \label{p1eq5}
(z-p)\left ((z-p)-\frac{1-p^2}{p}\right )\sum_{n=-1}^{\infty}a_n(z-p)^n \hspace{3cm}\\
= p+(z-p)-\frac{p}{1+p^2}\left (1+\sum_{n=0}^{\infty}c_n(z-p)^n
\right)\left((z-p)^2+2p(z-p)+p^2\right ).\nonumber
\eeq
Comparing  the coefficient of $(z-p)$  on both sides of (\ref{p1eq5}), we see that
\begin{equation}\label{p1eq6}
a_{-1}-\frac{1-p^2}{p}a_0 = \frac{1-p^2}{1+p^2}-\frac{p^2}{1+p^2}(2c_0+pc_1).
\end{equation}
Using the classical Schwarz-Pick lemma, it follows that
$$
|\omega'(p)|\leq \frac{1-|\omega (p)|^2}{1-p^2},\quad \mbox{i.e. }~|c_1|\leq \frac{1-|c_0|^2}{1-p^2}.
$$
In view of this observation, we have the estimate
$$
|2c_0+pc_1|\leq \frac{p(1-|c_0|^2)+2(1-p^2)|c_0|}{1-p^2}.
$$
For convenience, we set $x=|c_0|$ and consider
$$R_p(x)=p(1-x^2)+2(1-p^2)x.
$$
We see that $R_p(x)$ attains local maximum at $x_m=\frac{1-p^2}{p}$. Since
$x_m \geq 1$ for $p\in \left(0,\frac{\sqrt 5-1}{2}\right]$,
we see that
$$
|R_p(x)|\leq R_p(1)= 2(1-p^2)\quad, \, x\in[0,1],\quad p\in\left(0,\frac{\sqrt 5-1}{2}\right],
$$
and therefore, we have the estimate $|2c_0+pc_1|\leq 2$ for those $p$ in the said
interval. Now using this we get
from (\ref{p1eq6}) the  estimate (\ref{p1eq2}).
It is a simple exercise to see that the
equality is attained in (\ref{p1eq2}) for the following function
$$
f(z)=\frac{-zp}{(z-p)(1-zp)}.
$$
\hfill $\Box$

\begin{nonsec}Proof of Theorem \ref{1th2}\end{nonsec}
Let $f\in Co(p)$, with the expansion (\ref{p1eq1}).
Next, following the notation of the proof of Theorem \ref{1th1},
we compare the coefficients of $(z-p)^n$ $(n\geq 3$) on both side of
the equation (\ref{p1eq5}). This gives
$$a_{n-2}-\frac{1-p^2}{p}a_{n-1}=-\frac{p}{1+p^2}(c_{n-2}+2p c_{n-1}+p^2c_n)
\quad \mbox{$(n\geq 3)$}.
$$
Now, for a unimodular bounded analytic function $\omega$ in the unit disk $\ID$
having the expansion (\ref{p1eq4a}) in $\Delta_p$, we recall the following
result due to Ruscheweyh \cite[Theorem 2]{Rus-89}
$$
(1-p)^n(1+p)|c_n|\leq 1-|c_0|^2  \quad (n\geq 1)
$$
where the equality holds for $\omega(z)=e^{i\theta}, \theta\in[0,2\pi)$.
Using this, we easily obtain that
$$
\left|a_{n-2}-\frac{(1-p^2)a_{n-1}}{p}\right|
\leq \frac{p(1-|c_0|^2)}{(1+p^2)(1+p)(1-p)^n} \quad (n\geq 3).
$$
Consequently, (\ref{p1eq18}) follows since
$$
c_0=\frac{1-p^4}{p^4}a_{-1}+\frac{1}{p^2},
$$
by comparing  the constant terms on both sides of (\ref{p1eq5}).
Now, the equality holds in (\ref{p1eq18}) for the functions $f_\theta$,
$\theta\in [0, 2\pi)$,
in (\ref{p1eq18aa}), since both sides of the inequality are zero.
\hfill $\Box$

\section{Meromorphically starlike functions}\label{sec1-bpw3}
Let $\Sigma^s(p,w_0)$ denote the class of
meromorphic and univalent functions $f$ in $\ID$ (with the standard
normalization $f(0)=f'(0)-1=0$) having a simple pole
at $p\in (0,1)$ with the expansion (\ref{p1eq1}) such that $f$ is starlike
with respect to a fixed
$w_0\in \IC, w_0\neq 0$ (i.e.  $\overline\IC\backslash f(\ID)$ is a starlike
set with respect to $w_0$).
A well-known fact is that \cite{Living-94} $f\in\Sigma^s(p,w_0)$
if and only if ${\rm Re}\,\psi(z,f) >0$ for all $z\in \ID$, where
\be\label{p1eq20}
\psi(z,f) =\frac{-(z-p)(1-pz)f'(z)}{f(z)-w_0},\quad z\in \ID.
\ee

We now state and prove a useful representation formula for functions
in the class $\Sigma^s(p, w_0)$.

\bthm\label{1th5}
For $0<p<1$, let $f\in \Sigma^s(p, w_0)$. Then there exists a function
$\omega$ holomorphic in $\ID$ such that
$\omega(\ID)\subset \overline \ID,\, \omega(0)= -\frac{1}{2} (\frac{1}{w_0}+p+\frac{1}{p})$
and
\be\label{p1eq22}
f(z)= w_0 + \frac{pw_0(1+z\omega(z))^2}{(z-p)(1-zp)},\quad z\in \ID.
\ee
\ethm\bpf
The proof of this theorem is indeed a direct consequence of \cite[Corollary 2]{Zhang-92}
in which they use the notation $\sigma^*(p, w_0)$ in place of $\Sigma^s(p, w_0)$.
By this corollary we get if $f\in \Sigma^s(p, w_0)$, then
$$\left|\left\{\frac{f(z)-w_0}{pw_0}(z-p)(1-pz)\right\}^{1/2}-1\right|\leq |z| ,\quad z\in \ID.
$$
Now writing
$$
\omega(z)= \frac{1}{z}\left\{\frac{f(z)-w_0}{pw_0}(z-p)(1-pz)\right\}^{1/2}-\frac{1}{z}
$$
and simplifying the above expression for $f$ we get the desired
representation formula for functions in the class $\Sigma^s(p, w_0)$. Here we note
that $ \omega$ is holomorphic in $\ID$ and $| \omega(z)|\leq 1$. Also since $f'(0)=1$
we get $\omega(0)= -\frac{1}{2} (\frac{1}{w_0}+p+\frac{1}{p})$.
\epf

As a consequence of Theorem \ref{1th5}, we have the following result
which has been proved in \cite{Chang-88} but by using a different method of
proof.

\bcor
For $0<p<1$, let  $f\in \Sigma^s(p, w_0)$. Then, we have
\be\label{p1eq25}
\left |w_0+\frac{p(1+p^2)}{(1-p^2)^2}\right |\leq \frac{2p^2}{(1-p^2)^2} .
\ee
In particular, one has
$$\frac{p}{(1+p)^2}\leq |w_0|\leq \frac{p}{(1-p)^2}.
$$
\ecor \bpf
As $|\omega(0)|\leq 1$ and $\omega(0)= -\frac{1}{2} (\frac{1}{w_0}+p+\frac{1}{p})$,
it follows that
$$\left | \frac{1}{w_0}+\frac{1+p^2}{p}\right |\leq 2
$$
which is easily seen to be equivalent to the inequality (\ref{p1eq25}). The second inequality is a simple
consequence of  (\ref{p1eq25}).
\epf

\bthm\label{1th6}
Let $f\in\Sigma^s(p, w_0)$ have the Laurent expansion $(\ref{p1eq1})$. Then
\begin{enumerate}
\item[(i)] $\displaystyle
\left|a_{-1}-\frac{pw_0}{1-p^2}\right|\leq \frac{p|w_0|}{1-p^2}
\frac{|\frac{p}{w_0}+p^2+1|+2p^2}{2+|\frac{p}{w_0}+p^2+1|}
\left(\frac{|\frac{p}{w_0}+p^2+1|+2p^2}{2+|\frac{p}{w_0}+p^2+1|}
+2\right)$ ~, $p\in (0,1)$
\item[(ii)] $\displaystyle
 \left |a_0-\frac{1-p^2+p^4}{(1-p^2)^2}w_0\right|\leq \frac{p(2+2p-p^3)}{(1-p^2)^2}|w_0|
$~, $ p\in \left(0,\frac{\sqrt5 -1}{2}\right]$.
\end{enumerate}
Both inequalities are sharp for
$$
f(z)=\frac{-zp}{(z-p)(1-pz)}= w_0+\frac{pw_0}{(z-p)(1-pz)}(1-z)^2\in \Sigma^{s}(p, w_0)
$$
where $w_0=\frac{-p}{(1+p)^2}$.
\ethm\bpf
Let $f\in\Sigma^s(p, w_0)$ have the Laurent expansion (\ref{p1eq1})
and  consider the Taylor expansion  for $\omega$:
$$\omega(z)=\sum_{n=0}^{\infty}c_n(z-p)^n, \quad z\in\Delta_p.
$$
Now substituting (\ref{p1eq1}) and (\ref{p1eq4a}) in the representation formula
(\ref{p1eq22}) we get the following  series formulation of  (\ref{p1eq22}) valid in $\Delta_p$:\\

$\ds \sum_{n=-1}^{\infty} a_n(z-p)^n-w_0$
\beq\label{p1eq23}
&=& \frac{pw_0}{1-p^2}\sum_{n\geq 0}\left(\frac{p}{1-p^2}\right)^n (z-p)^{n-1}\\
\nonumber && \left[1+\{(z-p)^2+p^2+2p(z-p)\}\sum_{n\geq 0}\right.\left(\sum_{k=0}^{n}c_k c_{n-k}\right)(z-p)^n\\
\nonumber && \left.+\{2p+2(z-p)\}\sum_{n\geq 0}c_n(z-p)^n\right].
\eeq
While obtaining the above series form, we make use of the following relation
\beqq
(z-p)(1-pz) &=& (1-p^2)(z-p)\left (1- \frac{p}{1-p^2}(z-p) \right )\\
(1+z\omega(z))^2 &=& 1+ 2(z-p+p)w(z)+ ((z-p)^2+2p(z-p) +p^2)w(z)w(z).
\eeqq
Now, we proceed to prove (i).
Comparing the coefficients of $1/(z-p)$ on both sides of $(\ref{p1eq23})$ we get
$$a_{-1}=\frac{pw_0}{1-p^2}[1+p^2c_0^2+2pc_0].
$$
Finally, Schwarz Pick lemma applied to $\omega$ shows that
$$|c_0|=|\omega(p)|\leq \frac{|\omega(0)|+p}{1+p|\omega(0)|},
$$
where
$$|\omega(0)|=\frac{1}{2}\left|\frac{1}{w_0}+p+\frac{1}{p}\right| .
$$
Using this, we now get the desired estimate for $a_{-1}$.
It is also easy to check that the estimate stated in (i)
is sharp for the function mentioned in the statement of the theorem.\\

(ii) Comparing constant terms on  both sides of (\ref{p1eq23}) we get
$$a_0-w_0
= \frac{p^2 w_0}{(1-p^2)^2}(1+p^2c_0^2+2pc_0)+
\frac{2pw_0}{1-p^2}(p^2c_0c_1+p^2c_0^2+pc_1+c_0)
$$
or equivalently,
$$
a_0-\frac{1-p^2+p^4}{(1-p^2)^2}w_0
= \frac{p^2 w_0}{(1-p^2)^2}(p^2c_0^2+2pc_0)+
\frac{2pw_0}{1-p^2}(p^2c_0c_1+p^2c_0^2+pc_1+c_0).
$$
Now, we recall the well known estimates from the classical Schwarz Pick lemma:
$$|c_0|\leq 1, \quad |c_1|\leq \frac{1-|c_0|^2}{1-p^2}.
$$
For convenience, we use the notation $x= |c_0|$.
Using the above estimates, it is easy to see that the last equality
implies that
$$\left |a_0-\frac{1-p^2+p^4}{(1-p^2)^2}w_0\right|\leq
\frac{p|w_0|}{(1-p^2)^2}(2p+2x+2p^2x-2p^2x^3-p^3x^2).
$$
Next, we introduce
$$Q_p(x)= 2p+2x+2p^2x-2p^2x^3-p^3x^2, \quad 0\leq x\leq 1.
$$
It follows that $Q_p$ attains local maximum at
$$x_m= (-p^2+\sqrt{p^4+12(1+p^2)}\,)/(6p).
$$
Since $x_m\geq 1$ for
$p\in \left(0,\frac{\sqrt5 -1}{2}\right]$, we have
$$
\max\{Q_p(x):\, x\in[0,1]\}= Q_p(1)= 2+2p-p^3.
$$
This proves the inequality (ii) and the sharpness part can easily be
verified for the function given in the statement of
the theorem.
\epf

\br
It is a simple exercise to see that
\be\label{p1eq26}
\frac{\left |\frac{p}{w_0}+p^2+1\right |+2p^2}{2+\left |\frac{p}{w_0}+p^2+1
\right |}\leq p
\ee
is equivalent to
$$\frac{1}{2}\left |\frac{p}{w_0}+p^2+1\right |\leq p,
~\mbox{ i.e. }~ |\omega (0)|\leq 1 .
$$
Thus, (\ref{p1eq26}) holds.
If we use the  inequality (\ref{p1eq26}) then the inequality (i)
in Theorem \ref{1th6} turns out to be
$$
\left|a_{-1}-\frac{pw_0}{1-p^2}\right|\leq \frac{p^2}{1-p^2}(p+2)|w_0|,\quad p\in(0,1).
$$
Applying the triangle inequality in the above inequality and
the inequality (ii) in Theorem \ref{1th6} we get
\be\label{p1eq24}
|a_{-1}|\leq \frac{p(1+p)}{1-p}|w_0|, \quad p\in (0,1)
\ee
and
$$
|a_0|\leq \frac{1}{(1-p)^2}|w_0|, \quad p\in \left(0,\frac{\sqrt5 -1}{2}\right],
$$
respectively. Both the above estimates are sharp for the function stated in
Theorem \ref{1th6}.
\hfill $\Box$\er

In view of the estimate (\ref{p1eq24}) we observe that there was a minor error in
one of the results of Livingston, namely Theorem 9 in  \cite{Living-94}.
Indeed a counterexample is given by the function
$$g(z)=\frac{-zp}{(z-p)(1-pz)}\in \Sigma^{s}\left(p, \frac{-p}{1+p^2}\right).
$$
Here we note that
$$a_{-1}(g)= \frac{-p^2}{1-p^2}
$$
does not belong to the disk stated in Theorem 9 of \cite{Living-94}.
Moreover, the error is actually occurred in  \cite[p.~290]{Living-94} where the
inequality in the 6th line needs to be reversed, since $\xi-p\leq 0$.
We may now formulate a corrected version of \cite[Theorem 9]{Living-94}
for future use.

\bthm\label{1th7}
If $f\in\Sigma^{s}(p, w_0)$ and has the Laurent expansion
$(\ref{p1eq1})$, then we have
$$|a_{-1}|\geq \frac{p(1-p)}{1+p}|w_0|.
$$
The inequality is sharp for the function
$$g(z)=\frac{-zp}{(z-p)(1-pz)}
= w_0+\frac{pw_0}{(z-p)(1-pz)}(1-z)^2\in \Sigma^{s}(p, w_0)
$$
where $w_0=\frac{-p}{(1-p)^2}$.
\ethm

Here we also note that
$${\rm Re}\, \left(\frac{(z-p)(1-zp)g'(z)}{f(z)+\frac{p}{(1-p)^2}}\right)
=-(1-p)^2 {\rm Re}\left(\frac{1+z}{1-z}\right)<0
$$
for all $z\in \ID$ and $g$ satisfies the normalization condition $g(0)=0=g'(0)-1$
whenever $w_0=\frac{-p}{(1-p)^2}$.

\br
In view of the last theorem, the corollary that follows from Theorem 9 in \cite{Living-94}
is also not true since it uses the incorrect estimate for $|a_{-1}|$.
\hfill $\Box$
\er


{\bf Acknowledgement:} The authors thank Prof. K.-J. Wirths for his
suggestions on various stages of the paper and for his
continuous encouragement.


\end{document}